\documentclass[a4paper]{amsart}
\usepackage{amsmath,amssymb}
\usepackage[dvips]{graphics}
\usepackage[all]{xy}

\newtheorem{Theorem}{Theorem}[section]

\newtheorem{Proposition}[Theorem]{Proposition}

\makeatletter
\@addtoreset{figure}{section}%{subsection}
\def\@thmcountersep{-}
\makeatother

\numberwithin{equation}{section}

\begin{document}

\title[Unknotting numbers and crossing numbers of spatial planar graphs]{Unknotting numbers and crossing numbers of spatial embeddings of a planar graph}

\author{Yuta Akimoto}
\address{Graduate School of Education, Waseda University, Nishi-Waseda 1-6-1, Shinjuku-ku, Tokyo, 169-8050, Japan}
\email{motyanojika@gmail.com}

\author{Kouki Taniyama}
\address{Department of Mathematics, School of Education, Waseda University, Nishi-Waseda 1-6-1, Shinjuku-ku, Tokyo, 169-8050, Japan}
\email{taniyama@waseda.jp}

\thanks{The second author was partially supported by Grant-in-Aid for Scientific Research(A) (No. 16H02145) , Japan Society for the Promotion of Science.}


\date{}

\dedicatory{}


\begin{abstract}

It is known that the unknotting number $u(L)$ of a link $L$ is less than or equal to half the crossing number $c(L)$ of $L$. We show that there are a planar graph $G$ and its spatial embedding $f$ such that the unknotting number $u(f)$ of $f$ is greater than half the crossing number $c(f)$ of $f$. We study relations between unknotting number and crossing number of spatial embedding of a planar graph in general.
\end{abstract}

\maketitle

\section{Introduction}\label{introduction} 

Let $L$ be a link in the $3$-dimensional Euclidean space ${\mathbb R}^{3}$. Let $c(L)$ be the crossing number of $L$, that is, the minimal number of crossing points among all regular diagrams of $L$. Let $u(L)$ be the unknotting number of $L$, that is, the minimal number of crossing changes from $L$ to a trivial link. It is well-known that $u(L)$ is less than or equal to half of $c(L)$. See for example \cite{Taniyama2}. In this paper we show that this is not extended to spatial embeddings of planar graphs. We also show that this is extended to spatial embeddings of {\it trivializable} planar graphs. 

Let $G$ be a planar graph, that is, $G$ has an embedding into the plane ${\mathbb R}^{2}$, or equivalently, into the $2$-sphere ${\mathbb S}^{2}$. We denote the set of all vertices by ${\rm V}(G)$ and the set of all edges by ${\rm E}(G)$. 
An embedding $f:G\to{\mathbb R}^{3}$ is said to be a {\it spatial embedding} of $G$ and its image $f(G)$ is said to be a {\it spatial graph}. Let ${\rm SE}(G)$ be the set of all spatial embeddings of $G$. 
Let $\pi:{\mathbb R}^{3}\to{\mathbb R}^{2}$ be a natural projection defined by $\pi(x,y,z)=(x,y)$. 
We may assume up to ambient isotopy of ${\mathbb R}^{3}$ that $f\in{\rm SE}(G)$ is in general position with respect to $\pi$. 
Namely the composition map $\pi\circ f:G\to{\mathbb R}^{2}$ is a generic immersion. Here a continuous map $\varphi:G\to{\mathbb R}^{2}$ is said to be a {\it generic immersion} if it has only finitely many multiple points each of which is a transversal double point of an edge or two edges of $G$. Such a double point is said to be a {\it crossing point}, or simply a {\it crossing}. 
The set of all crossings of $\varphi$ is denoted by ${\mathcal C}(\varphi)$ and the number of crossings of $\varphi$ is denoted by $c(\varphi)=|{\mathcal C}(\varphi)|$. 
A {\it regular diagram}, or simply a {\it diagram} $D$ of $f\in{\rm SE}(G)$ is the image $\pi\circ f(G)\subset{\mathbb R}^{2}$ together with vertex/edge labels and over/under information at each crossing point. 
The set of all crossing points of $D$ is denoted by ${\mathcal C}(D)$ and the number of crossings of $D$ is denoted by $c(D)=|{\mathcal C}(D)|$. 

An element $t\in{\rm SE}(G)$ is said to be {\it trivial}, or {\it unknotted}, if there is a $2$-sphere $S$ embedded in ${\mathbb R}^{3}$ such that $t(G)\subset S$. It is known in \cite{Mason} that any two trivial embeddings of a planar graph are ambient isotopic. 
Therefore unknotting number is naturally extended to spatial embeddings of planar graphs as follows. 
For $f\in{\rm SE}(G)$, the {\it unknotting number} $u(f)$ of $f$ is defined to be the minimal number of crossing changes from $f$ to a trivial embedding of $G$. 
The {\it crossing number} $c(f)$ of $f$ is defined to be the minimal number of crossing points among all regular diagrams of spatial embeddings that are ambient isotopic to $f$. 

\begin{figure}[htbp]
      \begin{center}
\scalebox{0.5}{\includegraphics*{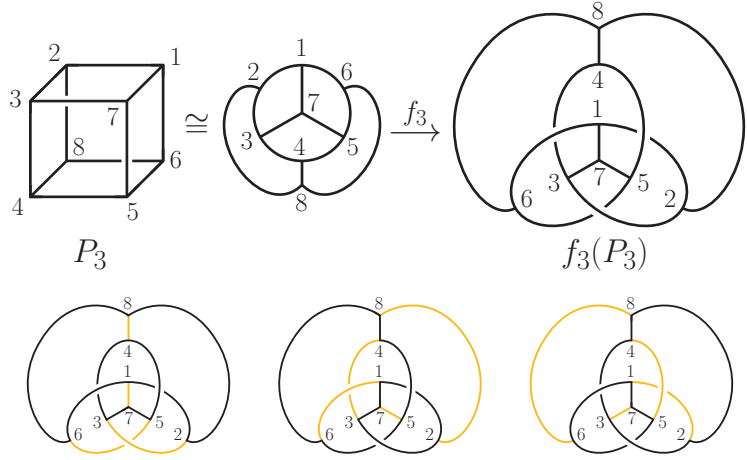}}
      \end{center}
   \caption{An example}
  \label{cube-graph}
\end{figure}

In comparison to link case stated above, we will show that there are a planar graph $G$ and a spatial embedding $f$ of $G$ such that $u(f)$ is greater than half of $c(f)$. 
Before giving an infinite sequence of such examples in Theorem \ref{theorem-plum}, we exhibit the first term of the sequence here. 
Let $P_{3}$ be the cube graph and $f_{3}\in {\rm SE}(P_{3})$ a spatial embedding of $P_{3}$ illustrated in Figure \ref{cube-graph}. 
Since $f_{3}(P_{3})$ contains a trefoil knot whose crossing number is $3$, we have $c(f_{3})=3$. The cube graph $P_{3}$ contains exactly $3$ pairs of mutually disjoint cycles. Each of them is a pair of parallel squares of the cube. Each of them forms a Hopf link in $f_{3}(P_{3})$ as illustrated in Figure \ref{cube-graph}. Therefore $f_{3}(P_{3})$ contains exactly $3$ Hopf links. Suppose that a crossing change is performed between two edges of $P_{3}$. Then only links that contain both of them may be changed. Since each edge of the cube is contained in exactly two squares, a crossing change may change at most $2$ Hopf links. 
Therefore at least one Hopf link remains after a crossing change and we have $u(f_{3})>1$. We will show later that $u(f_{3})\leq2$. 
Therefore $u(f_{3})=2$ and it is greater than half of $c(f_{3})=3$. 

\begin{figure}[htbp]
      \begin{center}
\scalebox{0.6}{\includegraphics*{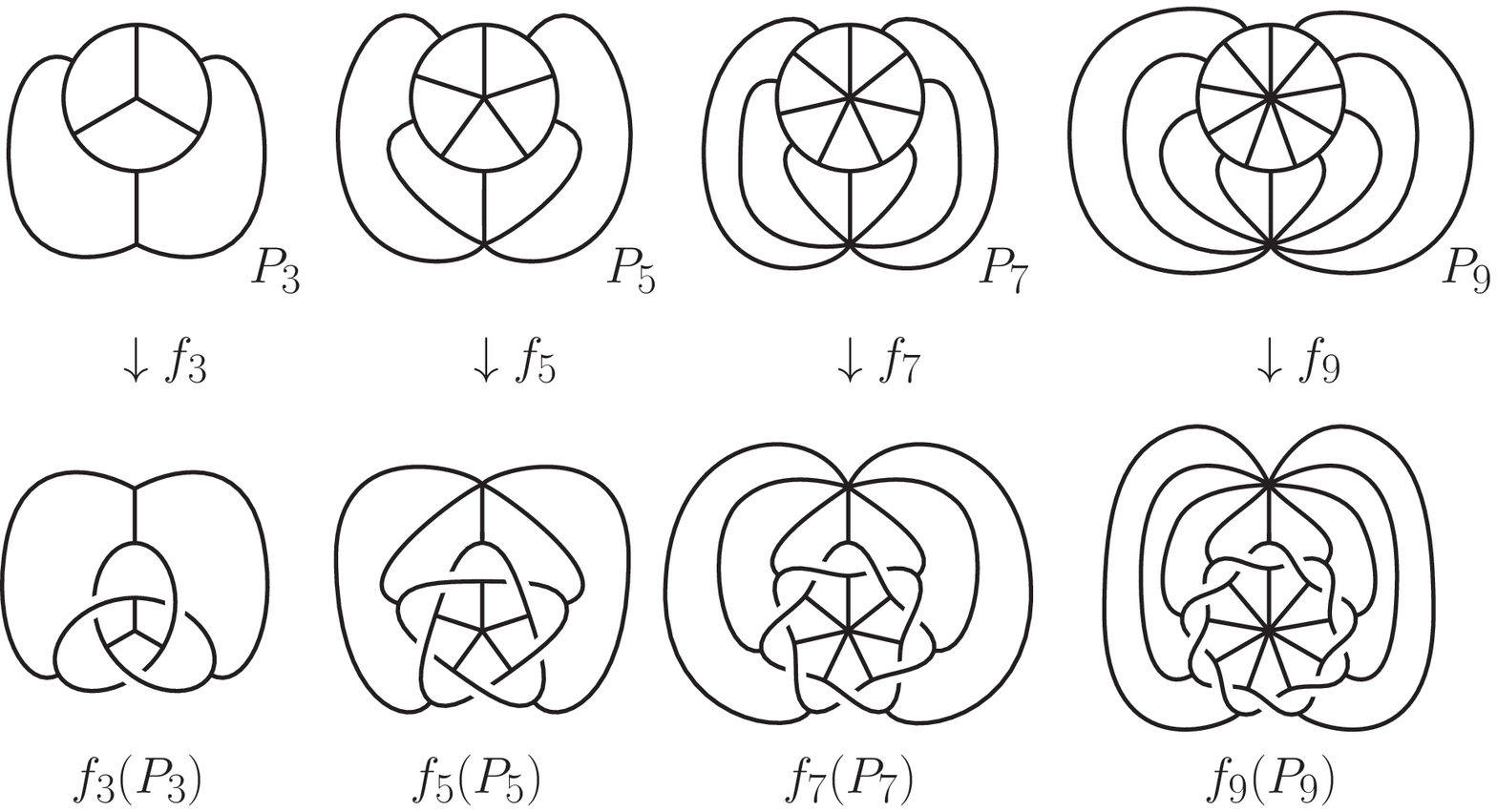}}
      \end{center}
   \caption{Examples}
  \label{plum-graphs}
\end{figure}

Let $n\in{\mathbb N}$ be a natural number. Let $P_{2n+1}$ be a planar graph and $f_{2n+1}\in {\rm SE}(P_{2n+1})$ a spatial embedding of $P_{2n+1}$ illustrated in Figure \ref{plum-graphs}. 
Since $f_{2n+1}(P_{2n+1})$ contains a $(2,2n+1)$-torus knot whose crossing number is $2n+1$, we have $c(f_{2n+1})=2n+1$. 
We have the following theorem whose proof will be given in \S \ref{proofs}.

\begin{Theorem}\label{theorem-plum}

Let $n\in{\mathbb N}$ be a natural number. Then $u(f_{2n+1})=2n$. 

\end{Theorem}

Now we explain the reason why it happens for these planar graphs. 
First we review the proof of $\displaystyle{u(L)\leq\frac{1}{2}c(L)}$ for a link $L$. 
It is well-known that every diagram of a link turns into a diagram of a trivial link by changing some of their crossings. Let $D$ be a diagram of $L$ with $c(D)=c(L)$. Then there is a subset $A$ of ${\mathcal C}(D)$ such that changing all crossings in $A$ turns $D$ to a diagram $T$ of a trivial link. Let $S$ be a diagram obtained from $D$ by changing all crossings in $B={\mathcal C}(D)\setminus A$. Note that $S$ is obtained from $T$ by changing all crossings of $T$. 
Therefore $S$ is a diagram of a mirror image of a trivial link. A mirror image of a trivial link is again a trivial link. 
Therefore $S$ is also a diagram of a trivial link. 
Thus we have $\displaystyle{u(L)\leq{\rm min}\{|A|,|B|\}\leq\frac{1}{2}c(D)=\frac{1}{2}c(L)}$. 

\begin{figure}[htbp]
      \begin{center}
\scalebox{0.5}{\includegraphics*{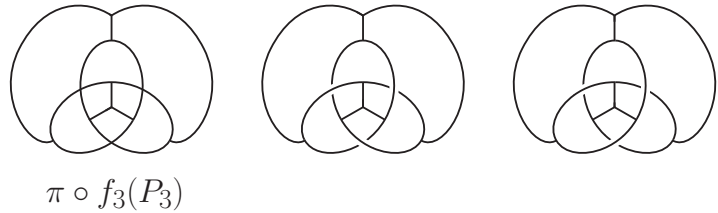}}
      \end{center}
   \caption{A knotted projection and diagrams coming from it}
  \label{knotted-projection}
\end{figure}

However there are a planar graph $G$ and a generic immersion $\varphi:G\to{\mathbb R}^{2}$ such that every spatial embedding $f$ of $G$ with $\pi\circ f=\varphi$ is not trivial \cite{Taniyama}. 
Such a generic immersion or its image is said to be a {\it knotted projection}. 
The projection $\pi\circ f_{3}(P_{3})$ illustrated in Figure \ref{knotted-projection} is a knotted projection first appeared in \cite{S-S}. 
Since $c(\pi\circ f_{3})=3$ there are $2^{3}=8$ diagrams coming from it. Up to mirror image and $2\pi/3$-rotation on the sphere ${\mathbb S}^{2}={\mathbb R}^{2}\cup\{\infty\}$ they are the diagrams illustrated in Figure \ref{knotted-projection}. One of them contains three Hopf links and the other contains one Hopf link.

A planar graph is said to be {\it trivializable} if it has no knotted projections. 
The set of all trivializable planar graphs is closed under minor reduction \cite{Taniyama}. 
A certain class of planar graphs are known to be trivializable \cite{Taniyama} \cite{S-S} \cite{Tamura} \cite{N-O-T-T}.

For a spatial embedding of a trivializable planar graph, the same argument as for a link works, and we have the following proposition. 

\begin{Proposition}\label{proposition-trivializable}

Let $G$ be a trivializable planar graph and $f:G\to{\mathbb R}^{3}$ a spatial embedding of $G$. Then $\displaystyle{u(f)\leq\frac{1}{2}c(f)}$.

\end{Proposition}

Then we are interested in the relation between unknotting number and crossing number of spatial embedding of a planar graph that holds in general.

\begin{Theorem}\label{theorem-linear-estimate}

Let $G$ be a planar graph. Then there exist real numbers $A$ and $B$ with the following property. 
For any spatial embedding $f:G\to {\mathbb R}^{3}$ of $G$, $\displaystyle{u(f)\leq A\cdot c(f)+B}$.

\end{Theorem}

\section{Proofs}\label{proofs} 
Some estimations below of the unknotting numbers of spatial embeddings of a planar graph are given in \cite{B-O}. 
In the following we give an estimation below based on an analysis of the changes of the linking numbers in a spatial graph. 

\vskip 5mm

\begin{figure}[htbp]
      \begin{center}
\scalebox{0.6}{\includegraphics*{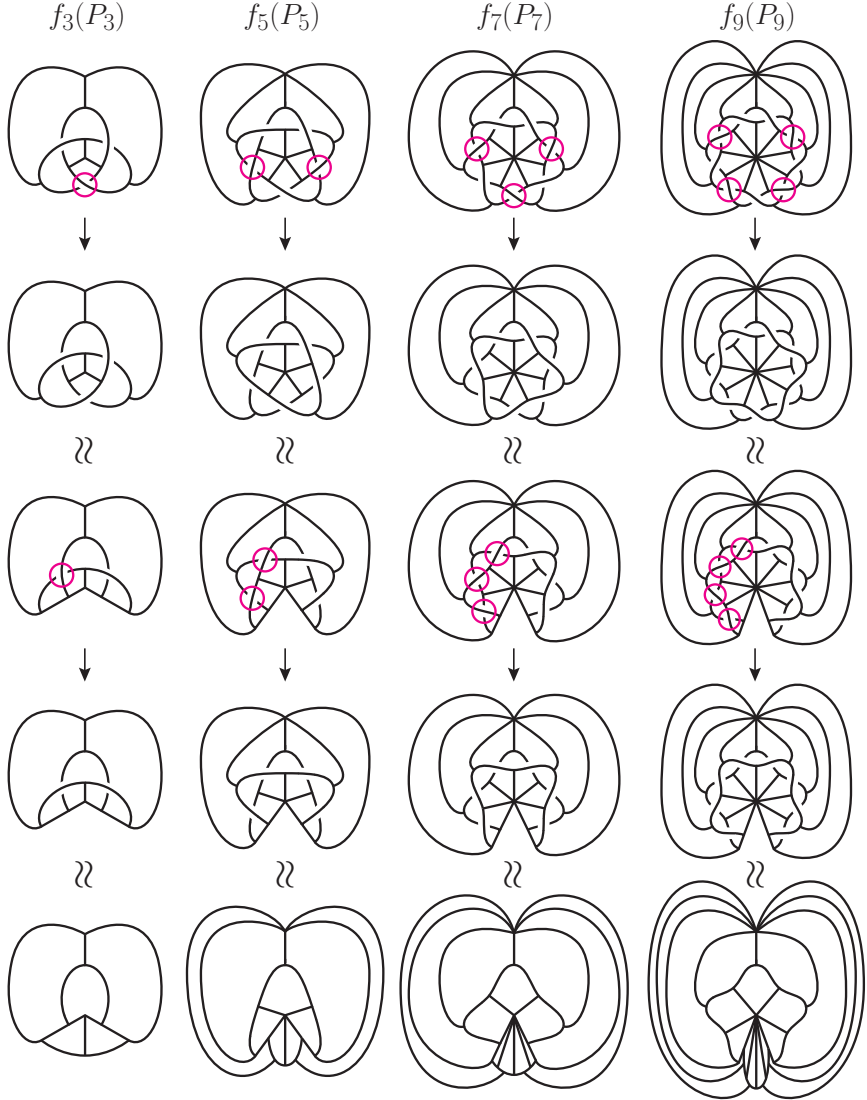}}
      \end{center}
   \caption{$u(f_{2n+1})\leq 2n$}
  \label{plum-unknotting}
\end{figure}

\noindent{\bf Proof of Theorem \ref{theorem-plum}.} First we show $u(f_{2n+1})\leq 2n$. It is shown in Figure \ref{plum-unknotting}.

\begin{figure}[htbp]
      \begin{center}
\scalebox{1.0}{\includegraphics*{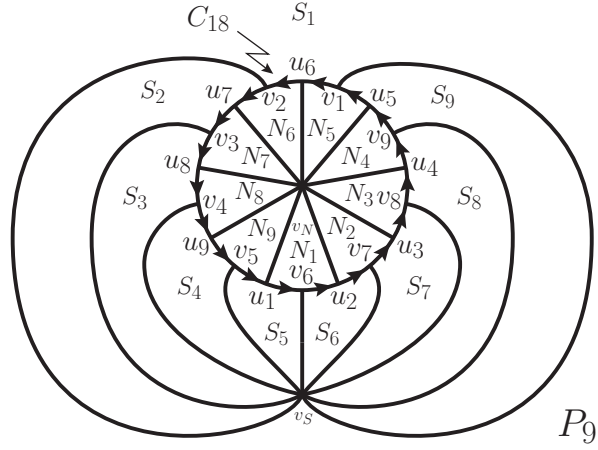}}
      \end{center}
   \caption{Cycles of $P_{9}$}
  \label{plum-cycles}
\end{figure}

Next we show $u(f_{2n+1})\geq 2n$. 
A prototype of the following proof is the proof of $u(f_{3})\geq2$ given in \S \ref{introduction}. 
Let $C_{4n+2}$ be the $4n+2$-cycle of $P_{2n+1}$ such that $f_{2n+1}(C_{4n+2})$ is a $(2,2n+1)$-torus knot. 
We call $C_{4n+2}$ the {\it equator} of $P_{2n+1}$. 
Each edge of $C_{4n+2}$ is said to be an {\it equatorial edge}. 
Let $v_{N}$ and $v_{S}$ the vertices of $P_{2n+1}$ not on $C_{4n+2}$. We call $v_{N}$ the {\it north pole} and $v_{S}$ the {\it south pole}. 
An edge incident to $v_{N}$ (resp. $v_{S}$) is said to be a {\it north spoke} (resp. {\it south spoke}). 
In the following we consider the suffixes of vertices of $P_{2n+1}$ modulo $2n+1$. 
Let $u_{1},u_{2},\cdots,u_{2n+1}$ be the vertices on $C_{4n+2}$ that are adjacent to $v_{N}$ and $v_{1},v_{2},\cdots,v_{2n+1}$ the vertices on $C_{4n+2}$ that are adjacent to $v_{S}$. 
We may assume that these vertices are cyclically arranged on $C_{4n+2}$ as $u_{1},v_{n+2},u_{2},v_{n+3},\cdots,u_{2n+1},v_{n+1}$. 
Let $N_{i}$ be a $4$-cycle $u_{i}v_{i+n+1}u_{i+1}v_{N}u_{i}$ oriented by this cyclic order of vertices and $S_{i}$ a $4$-cycle $v_{i}u_{i+n+1}v_{i+1}v_{S}v_{i}$ oriented by this cyclic order of vertices for $i=1,2,\cdots,2n+1$. 
These cycles are the region cycles when $P_{2n+1}$ is embedded into the $2$-sphere ${\mathbb S}^{2}$. 
Let ${\mathcal R}_{2n+1}=\{N_{1},\cdots,N_{2n+1},S_{1},\cdots,S_{2n+1}\}$ and we call ${\mathcal R}_{2n+1}$ the {\it region cycles} of $P_{2n+1}$. 
The case $n=4$ is illustrated in Figure \ref{plum-cycles}.

\noindent
For each spatial embedding $f$ of $P_{2n+1}$, we define the following integers. Here $\ell k(J,K)$ denotes the linking number of $J$ and $K$. 
\begin{align*}
\ell_{1}(f)&=\sum_{i=1}^{2n+1}\ell k(f(N_{i}),f(S_{i}))\\
\ell_{2}(f)&=\sum_{i=1}^{2n+1}\ell k(f(N_{i}),f(S_{i+1}))+\sum_{i=1}^{2n+1}\ell k(f(N_{i}),f(S_{i-1}))\\
\ell_{3}(f)&=\sum_{i=1}^{2n+1}\ell k(f(N_{i}),f(S_{i+2}))+\sum_{i=1}^{2n+1}\ell k(f(N_{i}),f(S_{i-2}))\\
\vdots&\\
\ell_{n}(f)&=\sum_{i=1}^{2n+1}\ell k(f(N_{i}),f(S_{i+n-1}))+\sum_{i=1}^{2n+1}\ell k(f(N_{i}),f(S_{i-n+1}))
\end{align*}
Then we define an element ${\mathcal L}(f)$ of ${\mathbb Z}^{n}$ by ${\mathcal L}(f)=(\ell_{1}(f),\ell_{2}(f),\cdots,\ell_{n}(f))$. 
We will observe how ${\mathcal L}(f)$ changes under a crossing change. 
Note that a crossing change between two strands of the same component does not change the linking number of a two-component link. Therefore a crossing change between two strands of the same edge does not change any linking numbers. Also a crossing change between two mutually adjacent edges does not change any linking numbers. Therefore we consider crossing changes between two mutually disjoint edges. 
Note that each edge of $P_{2n+1}$ is contained in exactly two region cycles of $P_{2n+1}$. Therefore a crossing change between two mutually disjoint edges of $P_{2n+1}$ changes at most four linking numbers. 
First we show an example. 
Suppose that a crossing change is performed between edge $v_{N}u_{1}$ and edge $v_{S}v_{1}$. 
Taking crossing signs under cycle orientations into account, we see that each of $\ell k(f(N_{1}),f(S_{1}))$ and $\ell k(f(N_{2n+1}),f(S_{2n+1}))$ changes by $\pm1$, each of $\ell k(f(N_{1}),f(S_{2n+1}))$ and $\ell k(f(N_{2n+1}),f(S_{1}))$ changes by $\mp1$, and all other linking numbers are unchanged. 
Therefore ${\mathcal L}(f)$ changes by $\pm(2,-2,0,\cdots,0)$. 
In such a way we have a list of possible changes of ${\mathcal L}(f)$ under a crossing change as follows. 

\begin{itemize}

\item Crossing changes between two equatorial edges
\begin{align*}
&\pm(2,0,\cdots,0),\pm(1,1,0,\cdots,0),\pm(0,2,0,\cdots,0),\pm(0,1,1,0,\cdots,0),\\
&\pm(0,0,2,0,\cdots,0),\pm(0,0,1,1,0,\cdots,0),\cdots\cdots,\\
&\pm(0,\cdots,0,2,0),\pm(0,\cdots,0,1,1),\pm(0,\cdots,0,2),\pm(0,\cdots,0,1)
\end{align*}

\item Crossing changes between an equatorial edge and a spoke
\begin{align*}
&\pm(1,-1,0,\cdots,0),\pm(0,1,-1,0,\cdots,0),\pm(0,0,1,-1,0,\cdots,0),\\
&\cdots\cdots,\pm(0,\cdots,0,1,-1,0),\pm(0,\cdots,0,1,-1),\pm(0,\cdots,0,1)
\end{align*}

\item Crossing changes between a north spoke and a south spoke
\begin{align*}
&\pm(2,-2,0,\cdots,0),\pm(1,-2,1,0,\cdots,0),\pm(0,1,-2,1,0,\cdots,0),\\
&\pm(0,0,1,-2,1,0,\cdots,0),\cdots\cdots,\pm(0,\cdots,0,1,-2,1,0),\\
&\pm(0,\cdots,0,1,-2,1),\pm(0,\cdots,0,1,-2),\pm(0,\cdots,0,1)
\end{align*}

\end{itemize}
Let ${\mathcal B}_{2n+1}$ be a finite subset of ${\mathbb Z}^{n}$ consists of the following elements. 
\begin{align*}
&a_{1}=(2,0,\cdots,0),a_{2}=(0,2,0,\cdots,0),\cdots,a_{n}=(0,\cdots,0,2)\\
&b_{1}=(1,1,0,\cdots,0),b_{2}=(0,1,1,0,\cdots,0),\cdots,b_{n-1}=(0,\cdots,0,1,1)\\
&c_{1}=(1,-1,0,\cdots,0),c_{2}=(0,1,-1,0,\cdots,0),\cdots,c_{n-1}=(0,\cdots,0,1,-1)\\
&d=e_{0}=(2,-2,0,\cdots,0)\\
&e_{1}=(1,-2,1,0,\cdots,0),e_{2}=(0,1,-2,1,0,\cdots,0),\cdots,e_{n-2}=(0,\cdots,0,1,-2,1)\\
&p=e_{n-1}=(0,\cdots,0,1,-2)\\
&q=b_{n}=c_{n}=e_{n}=(0,\cdots,0,1)
\end{align*}
We denote the $i$-th component of an element $v$ of ${\mathbb Z}^{n}$ by $v(i)$. Namely $v=(v(1),\cdots,v(n))$. 
Suppose that two spatial embeddings $f$ and $g$ of $P_{2n+1}$ are transformed into each other by $m$ times crossing changes. 
Then by the observation above we see that there is a map $\varphi:{\mathcal B}_{2n+1}\to{\mathbb Z}$ with 
$\displaystyle{
\sum_{v\in{\mathcal B}_{2n+1}}|\varphi(v)|\leq m
}$
such that
\[
{\mathcal L}(f)-{\mathcal L}(g)=\sum_{v\in{\mathcal B}_{2n+1}}\varphi(v)v.
\]
From Figure \ref{plum-graphs} we see that ${\mathcal L}(f_{2n+1})=(2n+1,0,\cdots,0)$. 
Note that ${\mathcal L}(t_{2n+1})=(0,\cdots,0)$ where $t_{2n+1}$ is a trivial embedding of $P_{2n+1}$. 
We have seen above that $f_{2n+1}$ and $t_{2n+1}$ are transformed into each other by $2n$ times crossing changes as illustrated in Figure \ref{plum-unknotting}. 
All of them are performed between two equatorial edges. 
First $n$ crossing changes correspond to $n$ times $(2,0,\cdots,0)$. Next $n$ crossing changes corresponds to 
\begin{align*}
&(1,1,0,\cdots,0),-(0,1,1,0,\cdots,0),(0,0,1,1,0,\cdots,0),-(0,0,0,1,1,0,\cdots,0),\cdots,\\
&(-1)^{n-3}(0,\cdots,0,1,1,0),(-1)^{n-2}(0,\cdots,0,1,1),(-1)^{n-1}(0,\cdots,0,1). 
\end{align*}
In fact the sum
\begin{align*}
&n(2,0,\cdots,0)+(1,1,0,\cdots,0)-(0,1,1,0,\cdots,0)+(0,0,1,1,0,\cdots,0)\\
&-(0,0,0,1,1,0,\cdots,0)+\cdots+(-1)^{n-3}(0,\cdots,0,1,1,0)\\
&+(-1)^{n-2}(0,\cdots,0,1,1)+(-1)^{n-1}(0,\cdots,0,1)
\end{align*}
is equal to $(2n+1,0,\cdots,0)={\mathcal L}(f_{2n+1})-{\mathcal L}(t_{2n+1})$. 

The following claim assures that at least $2n$ crossing changes are necessary. 

\vskip 5mm
\noindent
{\bf Claim.} {\it Let $\varphi:{\mathcal B}_{2n+1}\to{\mathbb Z}$ be a map with 
\[
\sum_{v\in{\mathcal B}_{2n+1}}\varphi(v)v=(2n+1,0,\cdots,0).
\]
Then 
$\displaystyle{
\sum_{v\in{\mathcal B}_{2n+1}}|\varphi(v)|\geq2n.
}$
}

\vskip 5mm
\noindent
We will show this claim step by step as follows. 

\vskip 5mm
\noindent
{\bf Subclaim 1.} {\it Let $\varphi:{\mathcal B}_{2n+1}\to{\mathbb Z}$ be a map with 
\[
\sum_{v\in{\mathcal B}_{2n+1}}\varphi(v)v=(2n+1,*,\cdots,*).
\]
Then 
$\displaystyle{
\sum_{v\in{\mathcal B}_{2n+1}}|\varphi(v)|\geq n+1.
}$
If 
$\displaystyle{
\sum_{v\in{\mathcal B}_{2n+1}}|\varphi(v)|=n+1,
}$
then $\displaystyle{\sum_{v\in{\mathcal B}_{2n+1}}\varphi(v)v(3)}$ is $0$ or $1$. 
}

\vskip 5mm
\noindent
{\bf Proof.} Note that $|v(1)|\leq2$ for any $v\in{\mathcal B}_{2n+1}$. Therefore $\displaystyle{
\sum_{v\in{\mathcal B}_{2n+1}}|\varphi(v)|\leq n
}$ implies 
$\displaystyle{
\sum_{v\in{\mathcal B}_{2n+1}}\varphi(v)v(1)\leq2n<2n+1.
}$
Thus we have 
$\displaystyle{
\sum_{v\in{\mathcal B}_{2n+1}}|\varphi(v)|\geq n+1. 
}$
Suppose 
$\displaystyle{
\sum_{v\in{\mathcal B}_{2n+1}}|\varphi(v)|=n+1. 
}$
Note that for $v\in{\mathcal B}_{2n+1}$, $v(1)=2$ if and only if $v\in\{a_{1},d\}$, $v(1)=1$ if and only if $v\in\{b_{1},c_{1},e_{1}\}$ and $v(1)=0$ otherwise. 
Then we see $\varphi(a_{1})+\varphi(d)=|\varphi(a_{1})|+|\varphi(d)|=n$ and $\varphi(b_{1})+\varphi(c_{1})+\varphi(e_{1})=|\varphi(b_{1})|+|\varphi(c_{1})|+|\varphi(e_{1})|=1$. 
Since $a_{1}(3)=d(3)=b_{1}(3)=c_{1}(3)=0$ and $e_{1}(3)=1$, we see that
\begin{align*}
\sum_{v\in{\mathcal B}_{2n+1}}\varphi(v)v(3)&=\varphi(a_{1})a_{1}(3)+\varphi(d)d(3)+\varphi(b_{1})b_{1}(3)+\varphi(c_{1})c_{1}(3)+\varphi(e_{1})e_{1}(3)
\end{align*}
is $0$ or $1$.
$\Box$

\vskip 5mm
\noindent
{\bf Subclaim 2.} {\it Let $\varphi:{\mathcal B}_{2n+1}\to{\mathbb Z}$ be a map with 
\[
\sum_{v\in{\mathcal B}_{2n+1}}\varphi(v)v=(2n+1,0,*,\cdots,*).
\]
Then 
$\displaystyle{
\sum_{v\in{\mathcal B}_{2n+1}}|\varphi(v)|\geq n+2.
}$
If 
$\displaystyle{
\sum_{v\in{\mathcal B}_{2n+1}}|\varphi(v)|=n+2,
}$
then $\displaystyle{\sum_{v\in{\mathcal B}_{2n+1}}\varphi(v)v(4)}$ is $0$, $1$ or $-1$. 
}

\vskip 5mm
\noindent
{\bf Proof.} By Subclaim 1 we have 
$\displaystyle{
\sum_{v\in{\mathcal B}_{2n+1}}|\varphi(v)|\geq n+1.
}$\\
Suppose 
$\displaystyle{
\sum_{v\in{\mathcal B}_{2n+1}}|\varphi(v)|=n+1.
}$
Then $\varphi(a_{1})+\varphi(d)=|\varphi(a_{1})|+|\varphi(d)|=n$ and $\varphi(b_{1})+\varphi(c_{1})+\varphi(e_{1})=|\varphi(b_{1})|+|\varphi(c_{1})|+|\varphi(e_{1})|=1$ as above. 
Since $a_{1}(2)=0,d(2)=-2,b_{1}(2)=1,c_{1}(2)=-1$ and $e_{1}(2)=-2$, we see that
\begin{align*}
\sum_{v\in{\mathcal B}_{2n+1}}\varphi(v)v(2)&=\varphi(a_{1})a_{1}(2)+\varphi(d)d(2)+\varphi(b_{1})b_{1}(2)+\varphi(c_{1})c_{1}(2)+\varphi(e_{1})e_{1}(2)
\end{align*}
is an odd number or a negative number, and cannot equal $0$. 
Therefore we have
$\displaystyle{
\sum_{v\in{\mathcal B}_{2n+1}}|\varphi(v)|\geq n+2.
}$
\\
Suppose
$\displaystyle{
\sum_{v\in{\mathcal B}_{2n+1}}|\varphi(v)|=n+2.
}$
Suppose
$\displaystyle{
\sum_{v\in{\mathcal B}_{2n+1},v(4)\neq0}|\varphi(v)|\geq2.
}$
Let $\psi:{\mathcal B}_{2n+1}\to{\mathbb Z}$ be a map defined by 
\[
\psi(v)=
\begin{cases}
\varphi(v)&(v(4)=0)\\
0&(v(4)\neq0).
\end{cases}
\]
Then we have
$\displaystyle{
\sum_{v\in{\mathcal B}_{2n+1}}|\psi(v)|\leq n. 
}$
Since $v(4)\neq0$ implies $v(1)=0$, we have 
\[
\sum_{v\in{\mathcal B}_{2n+1}}\psi(v)v=(2n+1,*,\cdots,*).
\]
This contradicts Subclaim 1. 
Therefore 
$\displaystyle{
\sum_{v\in{\mathcal B}_{2n+1},v(4)\neq0}|\varphi(v)|\leq1.
}$
Suppose that $u(4)\neq0$ and $|\varphi(u)|=1$. 
Suppose that $u\neq e_{2}=(0,1,-2,1,0,\cdots,0)$. 
Then $u(1)=u(2)=0$. 
Let $\psi:{\mathcal B}_{2n+1}\to{\mathbb Z}$ be a map defined by 
\[
\psi(v)=
\begin{cases}
\varphi(v)&(v\neq u)\\
0&(v=u).
\end{cases}
\]
Then we have
$\displaystyle{
\sum_{v\in{\mathcal B}_{2n+1}}|\psi(v)|=n+1 
}$
and 
\[
\sum_{v\in{\mathcal B}_{2n+1}}\psi(v)v=(2n+1,0,*,\cdots,*).
\]
This is a contradiction. Therefore $u=(0,1,-2,1,0,\cdots,0)$ and 
$\displaystyle{\sum_{v\in{\mathcal B}_{2n+1}}\varphi(v)v(4)}=\pm1$. 
$\Box$

\vskip 5mm
Now we inductively show the subclaims for $3\leq k\leq n$. 

\vskip 5mm
\noindent
{\bf Subclaim $k$.} {\it Let $\varphi:{\mathcal B}_{2n+1}\to{\mathbb Z}$ be a map with 
\[
\sum_{v\in{\mathcal B}_{2n+1}}\varphi(v)v=(2n+1,\underbrace{0,\cdots,0}_{k-1},*,\cdots,*).
\]
Then 
$\displaystyle{
\sum_{v\in{\mathcal B}_{2n+1}}|\varphi(v)|\geq n+k.
}$\\
If $k\leq n-2$ and 
$\displaystyle{
\sum_{v\in{\mathcal B}_{2n+1}}|\varphi(v)|=n+k,
}$
then $\displaystyle{\sum_{v\in{\mathcal B}_{2n+1}}\varphi(v)v(k+2)}$ is $0$, $1$ or $-1$. 
}

\vskip 5mm
\noindent
{\bf Proof.} By Subclaim $k-1$ we have 
$\displaystyle{
\sum_{v\in{\mathcal B}_{2n+1}}|\varphi(v)|\geq n+k-1.
}$\\
Suppose 
$\displaystyle{
\sum_{v\in{\mathcal B}_{2n+1}}|\varphi(v)|=n+k-1.
}$
Note that the sum 
\[
\sum_{i=1}^{k}\left(\sum_{v\in{\mathcal B}_{2n+1}}\varphi(v)v(i)\right)=(2n+1)+\underbrace{0+\cdots+0}_{k-1}=2n+1
\]
is odd. For $v\in{\mathcal B}_{2n+1}$, 
$\displaystyle{
\sum_{i=1}^{k}v(i)
}$
is odd if and only if $v\in\{b_{k},c_{k},e_{k-1},e_{k}\}$. 
Therefore there exist $u\in\{b_{k},c_{k},e_{k-1},e_{k}\}$ with $\varphi(u)\neq0$. 
Let $\psi:{\mathcal B}_{2n+1}\to{\mathbb Z}$ be a map defined by 
\[
\psi(v)=
\begin{cases}
\varphi(v)&(v\neq u)\\
0&(v=u).
\end{cases}
\]
Then we have
$\displaystyle{
\sum_{v\in{\mathcal B}_{2n+1}}|\psi(v)|\leq n+k-2. 
}$
Suppose $u\in\{b_{k},c_{k},e_{k}\}$. Then we see $u(i)=0$ for all $i\leq k-1$. 
Therefore we have 
\[
\sum_{v\in{\mathcal B}_{2n+1}}\psi(v)v=(2n+1,\underbrace{0,\cdots,0}_{k-2},*,\cdots,*).
\]
This contradicts Subclaim $k-1$. 
Suppose $u=e_{k-1}$. Then we have 
\[
\sum_{v\in{\mathcal B}_{2n+1}}\psi(v)v=(2n+1,\underbrace{0,\cdots,0}_{k-3},-\varphi(u),-2\varphi(u),*,\cdots,*).
\]
This contradicts Subclaim $k-2$ that asserts 
$\displaystyle{\sum_{v\in{\mathcal B}_{2n+1}}\psi(v)v(k)}$ is $0$, $1$ or $-1$. 
Therefore we have
$\displaystyle{
\sum_{v\in{\mathcal B}_{2n+1}}|\varphi(v)|\geq n+k.
}$
\\
Suppose 
$\displaystyle{
\sum_{v\in{\mathcal B}_{2n+1}}|\varphi(v)|=n+k.
}$
Suppose
$\displaystyle{
\sum_{v\in{\mathcal B}_{2n+1},v(k+2)\neq0}|\varphi(v)|\geq2.
}$
Let $\psi:{\mathcal B}_{2n+1}\to{\mathbb Z}$ be a map defined by 
\[
\psi(v)=
\begin{cases}
\varphi(v)&(v(k+2)=0)\\
0&(v(k+2)\neq0).
\end{cases}
\]
Then we have
$\displaystyle{
\sum_{v\in{\mathcal B}_{2n+1}}|\psi(v)|\leq n+k-2. 
}$
Since $v(k+2)\neq0$ implies $v(i)=0$ for all $i\leq k-1$, we have 
\[
\sum_{v\in{\mathcal B}_{2n+1}}\psi(v)v=(2n+1,\underbrace{0,\cdots,0}_{k-2},*,\cdots,*).
\]
This contradicts Subclaim $k-1$. 
Therefore 
$\displaystyle{
\sum_{v\in{\mathcal B}_{2n+1},v(k+2)\neq0}|\varphi(v)|\leq1.
}$
Suppose that $u(k+2)\neq0$ and $|\varphi(u)|=1$. 
Suppose that $u\neq e_{k}=(\underbrace{0,\cdots,0}_{k-1},1,-2,1,0,\cdots,0)$. 
Then $u(i)=0$ for all $i\leq k$. 
Let $\psi:{\mathcal B}_{2n+1}\to{\mathbb Z}$ be a map defined by 
\[
\psi(v)=
\begin{cases}
\varphi(v)&(v\neq u)\\
0&(v=u).
\end{cases}
\]
Then we have
$\displaystyle{
\sum_{v\in{\mathcal B}_{2n+1}}|\psi(v)|=n+k-1 
}$
and 
\[
\sum_{v\in{\mathcal B}_{2n+1}}\psi(v)v=(2n+1,\underbrace{0,\cdots,0}_{k-1},*,\cdots,*).
\]
This is a contradiction. Therefore we have $u=(\underbrace{0,\cdots,0}_{k-1},1,-2,1,0,\cdots,0)$ and 
$\displaystyle{\sum_{v\in{\mathcal B}_{2n+1}}\varphi(v)v(k+2)}=\pm1$. 
$\Box$

\vskip 5mm
Note that Subclaim $n$ is equal to Claim. Thus we have shown $u(f_{2n+1})\geq 2n$. 
$\Box$

\vskip 3mm

\begin{figure}[htbp]
      \begin{center}
\scalebox{0.6}{\includegraphics*{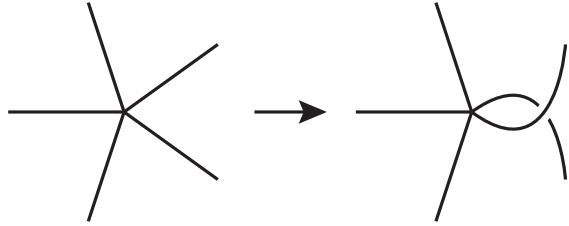}}
      \end{center}
   \caption{A deformation}
  \label{cyclic-order}
\end{figure}

\noindent{\bf Proof of Theorem \ref{theorem-linear-estimate}.} We choose and fix a spanning tree $T$ of $G$ and a vertex $v$ of $T$. 
We also choose a fixed embedding $G\subset{\mathbb R}^{2}$. Then we have an embedding of $T$ into ${\mathbb R}^{2}$ by $T\subset G\subset{\mathbb R}^{2}$. 
For a vertex $u$ of $G$, the degree of $u$ in $G$ is denoted by ${\rm deg}(u,G)$. For a vertex $u$ of $T$, the degree of $u$ in $T$ is denoted by ${\rm deg}(u,T)$. Note that ${\rm deg}(u,T)\leq{\rm deg}(u,G)$. 
Let $D$ be a diagram of $f$ with $c(D)=c(f)$. 
For each vertex $u$ of $T$ with ${\rm deg}(u,T)\geq3$ we repeatedly perform the deformation of $D$ near $u$ illustrated in Figure \ref{cyclic-order} such that the cyclic order of the edges of $T$ incident to $u$ near $u$ coincides with that of $u$ in $T\subset G\subset{\mathbb R}^{2}$. Note that the deformation involves not only the edges of $T$ but all edges of $G$ incident to $u$. 
Each deformation increase the number of crossings by $1$. 
We will see that the increase of $c(D)$ is bounded above by a constant $b$ that depends only on $G$. 
Suppose that ${\rm deg}(u,G)$ is $2k_{u}+1$ or $2k_{u}+2$. Let $l_{u}={\rm deg}(u,T)$. 
We fix an edge of $T$ incident to $u$ and deform other edges of $T$ incident to $u$ so that the cyclic order of them is equal to that in $T\subset G\subset{\mathbb R}^{2}$. 
By the choice of turning right or left, at most $k_{u}$ times application of the deformation illustrated in Figure \ref{cyclic-order} is sufficient to move each edge into a designated position. 
Therefore we see that at most $k_{u}(l_{u}-1)$ times application of the deformation illustrated in Figure \ref{cyclic-order} is sufficient to change the cyclic order of the edges of $T$ incident to $u$. 
Therefore at most $b=\displaystyle{\sum_{u}k_{u}(l_{u}-1)}$ times application is sufficient. 
Thus the diagram $D$ is deformed into a diagram $D'$ with $c(D')\leq c(D)+b$.

\begin{figure}[htbp]
      \begin{center}
\scalebox{0.6}{\includegraphics*{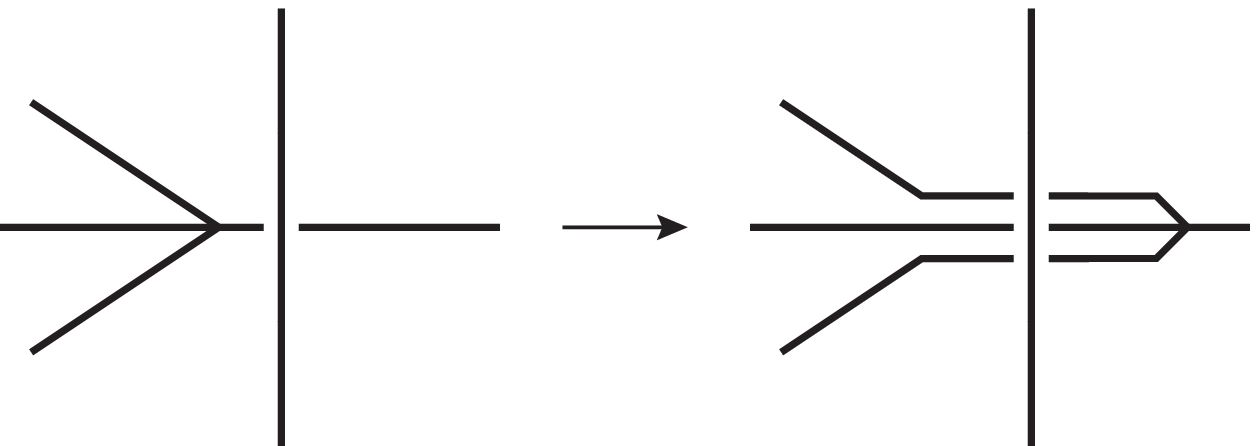}}
      \end{center}
   \caption{A deformation}
  \label{shrinking}
\end{figure}

\begin{figure}[htbp]
      \begin{center}
\scalebox{1.0}{\includegraphics*{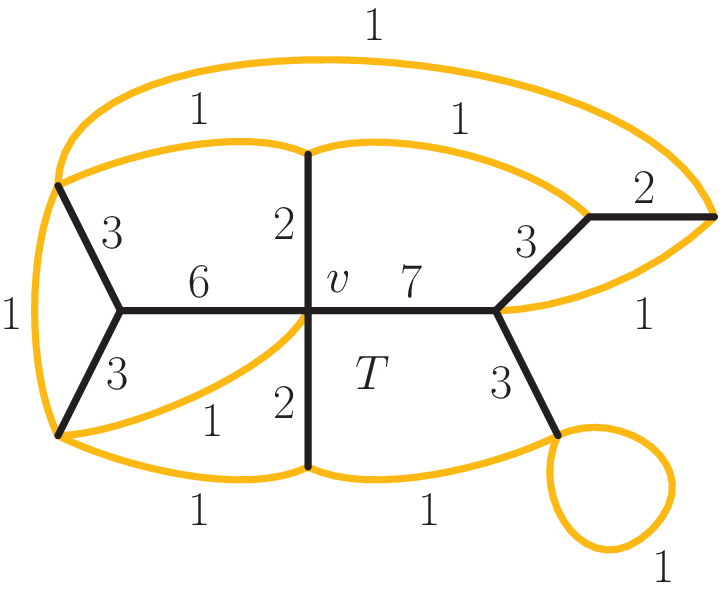}}
      \end{center}
   \caption{An example}
  \label{branch-index}
\end{figure}

Next we deform $D'$ to $D''$ as described below such that there are no crossings of $D''$ on $T$. 
We shrink the edges of $T$ toward $v$. 
In the process the number of crossings may increase as illustrated in Figure \ref{shrinking}. 
In Figure \ref{shrinking} only one case of over/under crossing information is illustrated. Another case is similar. 
To estimate the increase of the number of crossings, we define the following. 
Let $e$ be an edge of $T$. We denote the shortest path of $T$ containing $e$ and $v$ by $P_{T}(e,v)$. 
We define a positive integer $b(e)$ for each edge $e$ of $G$ as follows. Set $b(e)=1$ for $e\in E(G)\setminus E(T)$. 
Let $e$ be an edge of $T$. Let $u$ be a terminal vertex of $P_{T}(e,v)$ with $u\neq v$. Then $u$ is incident to $e$. Suppose that $b(d)$ is already defined for all edges but $e$ incident to $u$. Then $b(e)$ is defined to be the sum of these $b(d)$s where a loop is counted twice. This recursively defines $b(e)$ for all $e\in E(G)$. See for example Figure \ref{branch-index}. 
We see from Figure \ref{shrinking} that a crossing in $D'$ between edges $d$ and $e$ corresponds $b(d)b(e)$ crossings in $D''$. 
Let $a=({\rm max}\{b(e)|e\in E(G)\})^{2}$. Then we have $c(D'')\leq a\cdot c(D')$. Since $c(D')\leq c(D)+b$ we have $c(D'')\leq a\cdot c(D)+ab$. Note that $D''$ cannot be a knotted projection. In fact a descending algorithm under any ordering and orientation of the edges $E(G)\setminus E(T)$ produce a trivial embedding of $G$. Therefore we have $\displaystyle{u(f)\leq\frac{1}{2}c(D'')}$. Thus we have $\displaystyle{u(f)\leq\frac{1}{2}(a\cdot c(D)+ab)}$. Set $\displaystyle{A=\frac{1}{2}a}$ and $\displaystyle{B=\frac{1}{2}ab}$ we have $\displaystyle{u(f)\leq A\cdot c(f)+B}$. 
$\Box$


{\normalsize
\end{document}